\documentclass[10pt]{article}
\usepackage{latexsym} 
\usepackage{amsmath}
\usepackage{amssymb}

\usepackage[applemac ]{inputenc}
\usepackage[T1]{fontenc}

\usepackage{amsfonts }
\usepackage{amsthm}
\usepackage[mathcal]{eucal}
\usepackage{color}

\title{Stochastic approximations of set-valued dynamical systems: Convergence with positive probability to an attractor \thanks{We acknowledge financial support from the  Swiss National Science Foundation Grant 200021-103625/1}}

\author{
Mathieu Faure
\\{\small mathieu.faure@unine.ch}
\\ Gregory Roth 
\\ {\small gregory.roth@unine.ch}
\zdeux
\\{\small Institut de Math\'ematiques, Universit\'e de Neuch\^atel}, 
\\{\small Rue Emile-Argand 11. Neuch\^atel. Switzerland}.
}

\def\R{\mathbb R}

\def\hop{{\noindent}}

\newcommand{\vs}[1]{\vspace{#1}}

\newcommand{\begitem}{\begin{itemize}}
\newcommand{\finit}{\end{itemize}}

\newcommand{\zun}{\vs{0.1cm}} 
 
\newcommand{\zdeux}{\vs{0.2cm}} 
\newcommand{\ztrois}{\vs{0.3cm}}
\newcommand{\zcinq}{\vs{0.5cm}} 

\newcommand{\F}{\mathcal{F}}

\begin{document}

\thispagestyle{empty}

\newtheorem{theoreme}{Theorem}[section]
\newtheorem{lemme}[theoreme]{Lemma}
\newtheorem{proposition}[theoreme]{Proposition}
\newtheorem{hypothese}[theoreme]{Hypothesis}
\newtheorem{definition}[theoreme]{Definition}
\newtheorem{Remarque}[theoreme]{Remark}
\newtheorem{corollaire}[theoreme]{Corollary}
\newtheorem{exemple}[theoreme]{Example}
\newtheorem{axiome}[theoreme]{Axiom}
\newtheorem{condition}[theoreme]{Condition}
\newtheorem{example}[theoreme]{Example}
\def\bex{\begin{example} {\rm }
\def\eex{\end{example} }}

\maketitle

\begin{abstract}
\hop A successful method to describe the asymptotic behavior of a discrete time stochastic process governed by some recursive formula is to relate it to the limit sets of a well chosen mean differential equation. Under an attainability condition, convergence to a given attractor of the flow induced by this  dynamical system was proved to occur with positive probability (Bena\"im, 1999) for a class of Robbins Monro algorithms. Bena\"im et al. (2005) generalised this approach for stochastic approximation algorithms whose average behavior is related to a differential inclusion instead. We pursue the analogy by extending to this setting the result of convergence with positive probability to an attractor.
\end{abstract}
\zdeux

\hop {\bfseries Key words:} Stochastic approximations, set-valued dynamical systems, attractor, game theory, Markovian fictitious play
\zun

\hop {\bfseries MSC2000 Subject classification:} 62L20, 34A60, 34B40, 34B41, 91A25, 91A26

\hop {\bfseries OR/MS subject classification:} stochastic model applications

\section{Introduction}
\zcinq

\subsection{Settings and bibliography}
\ztrois

\hop Stochastic approximation algorithms were born in the early 50s through the work of Robbins and Monro \cite{RobMon51} and Kiefer and Wolfowitz \cite{KieWol52}. Consider a discrete time stochastic process $(x_n)_{n\ge0}$ defined  by the following recursive formula:
\begin{equation}\label{algo1}
x_{n+1} - x_n = \gamma_{n+1} \left(F(x_n) + U_{n+1}\right),
\end{equation}
where $F : \R^m \rightarrow \R^m$ is a Lipschitz function, $(\gamma_n)_{n}$ is a positive decreasing sequence  and $(U_n)_{n}$ a sequence of $\R^m$-valued random variables defined on a probability space $(\Omega, \F,P)$, adapted to some filtration $(\mathcal{F}_n)_n \subset \mathcal{F}$. In order to describe the limit behavior of the sample paths $(x_n(\omega))_n$, a natural idea is to compare them to the solution curves of the dynamical system induced by the ordinary differential equation
\begin{equation}\label{ode}
	\frac{dx}{dt} = F(x).
\end{equation}
This is the celebrated method of ordinary differential equation (ODE) which was introduced by Ljung in \cite{Lju77}. Heuristically, one can think of (\ref{algo1}) as a kind of Cauchy-Euler approximation scheme for numerically solving (\ref{ode}) with step size $(\gamma_n)_n$ and an added noise $(U_n)_{n}$. We could reasonably expect that, under appropriate assumptions on  $(\gamma_n)_n$ and if the noise $(U_n)_n$ vanishes, the asymptotic behaviors of $(x_n)_n$ and  the ODE are closely related.

\hop Thereafter, the method was studied and developed by many people (see Kushner and Clark \cite{KusCla78}, Benveniste et al \cite{BMP90},  Duflo \cite{Duf96} or Kushner and Yin \cite{KusYin03}). Originally, only simple dynamics were considered, for example the negative of the gradient of a cost function. However, it appears in several situations, for example, learning models or game theory, that the corresponding vector field may be more complex. 
\zdeux

\hop Bena\"im and Hirsch  have conducted, in a series of papers (essentially \cite{BenHir96} and \cite{Ben96}),  a thorough study of this method. They proved that the asymptotic behavior of stochastic approximation process can be described with a great deal of generality through the study of  the asymptotics of the ODE . One of the main results is the characterization of  limit sets $(x_n(\omega))_n$ via the flow induced by $ F $, in the sense that, almost surely, these sets are compact, invariant  and contain no proper attractor for the deterministic flow (this is the notion of \emph{internal chain recurrence} in the sense of Conley \cite{Con78}, see also Bowen \cite{Bow75}). 

\hop Now, let $F : \R^m \rightrightarrows \R^m$ be a sufficiently regular set-valued map and consider some discrete time stochastic processes $(x_n)_{n\ge0}$ satisfying the following recursive formula:
\begin{equation}\label{algoINC}
x_{n+1} - x_n - \gamma_{n+1} U_{n+1} \in \gamma_{n+1}F(x_n),
\end{equation}
where $(\gamma_n)_{n}$ is a positive decreasing sequence  and $(U_n)_{n}$ a sequence of $\R^m$-valued random variables defined on a probability space $(\Omega, \F,P)$.

\hop In \cite{BHS1}, Bena\"im, Hofbauer and Sorin have generalized  the ODE method to the algorithms given by (\ref{algoINC}) and extended  the characterization of limits set in the sense that they are again,  under certain assumptions on the step size and the noise, connected and attractor free
for the set-valued dynamical system induced by the differential inclusion

\begin{equation} \label{inclusiondiff}
\frac{dx}{dt} \in F(x).
\end{equation}

\hop This generalization allows us to extend this technique to a much wider class of problems arising, for example, in economics or game theory (see Bena\"im, Hofbauer and Sorin \cite{BHS2}). 
\zdeux

\hop In this paper, we pursue the analogy between the ODE method and the differential inclusion method. The aim is to extend to the case of differential inclusions, the result of Bena\"im (see theorem 7.3 in \cite{Ben99}) which guarantees that, under certain assumptions on the step size and the noise, the stochastic approximation process  converges with positive probability to a given attractor of the set-valued dynamical system induced by $F$.

\zdeux
\hop The organization of the paper is as follows. In section 1.2, we define a standard set-valued map and introduce the crucial notion of attainability so as to state a simple version of the main result. In section 2, we introduce the different notions of internal chain transitivity, asymptotic pseudotrajectories and perturbed solutions. Our main assumption (hypothesis \ref{H1}) is given, the convergence result is stated in full generality and we define a generalised stochastic approximation process which satisfies the above assumption. An example of adaptive learning process to which our results may be applied  is given in section 3. Finally, the proof of a crucial result needed in our study is postponed to section 4.

\subsection{The main result, a simple version}
\ztrois

\hop In the following, $M \subset \R^m$ is a compact set.
\begin{definition}\label{standarsetvaluedmap}[Standard set-valued map] A correspondance $F : \R^m \rightrightarrows \R^m$ is said to be \emph{standard} if it satisfies the following assumptions:

\begitem

\item for any $x \in \R^m$, $F(x)$ is a non empty, compact and convex set of $\R^m$,

\item $F$ is \emph{closed}, which means that its graph
\[Gr(F) := \left\{ (x,y) \in \R^m \times \R^m \mid \; \, y \in F(x)\right\}\]

\hop is closed,

\item there exists $c>0$ such that
\[\sup_{z \in F(x)} \|z\| \leq c (1 + \|x\|).\]

\end{itemize}

\end{definition}

\hop Under the above assumptions, it is well known (see Aubin and Cellina \cite{AubCel84}) that (\ref{inclusiondiff}) admits at least one solution (i.e. an absolutely continuous mapping $\mathbf{x} : \R \rightarrow \R^m$ such that $\dot{\mathbf{x}}(t) \in F(\mathbf{x}(t))$ for almost every $t$) through any initial point. 

\hop  We call $S_x$ the set of solutions with initial condition $\mathbf{x}(0) = x$. The set-valued dynamical system induced by the differential inclusion will be denoted $\Phi = (\Phi_t)_{t \in \mathbb{R}}$. To any $x \in \mathbb{R}^m$, it associates the non empty set
\[ \Phi_t(x) := \left\{ \mathbf{x}(t) \mid \; \mathbf{x} \in S_x \right\}.\]
Finally, $S_{\Phi} := \cup_{x} S_x$ is the set of all solution curves. In order to understand the main result, recall some classical definitions about the set-valued dynamics.

\begin{definition}
	A non empty compact set $A \subset \R^m$ is called an \emph{attractor} for $\Phi$, provided it is \emph{invariant} (i.e. for all $x\in A$, there exists a solution $\mathbf{x}$ to (\ref{inclusiondiff}) with $\mathbf{x}(0) =x$ and such that $\mathbf{x}(\R) \subset A$) and that there is a 		neighborhood $U$ of $A$ with the property that, for every $\epsilon >0$, there exists $t_{\epsilon}>0	$ such that
	$$\Phi_t(U) \subset N^{\epsilon}(A)$$
	for all $t\ge t_{\epsilon}$, where $N^{\varepsilon}(A)$ is the $\varepsilon$-neighborhood of $A$. An open set $U$ with this property is called a \emph{fundamental neighborhood} of $A$.
\end{definition}

\begin{definition}
Let $A \subset \R^m$ be an attractor for the set-valued dynamical system. The \emph{basin of attraction} of $A$ is the set $$\mathcal{B}(A) := \{x \in \R^m :  \Phi_{[0,+ \infty[}(x) \, \mbox{ bounded and } \; \omega_{\Phi}(x) \subset A \},$$ where $\omega_{\Phi}(x) =\bigcap_{t\ge 0}\overline{\Phi_{[t, \infty[}(x)}$ is the omega limit set of the point $x$. 
\end{definition}

\hop Now consider a discrete time stochastic process $(x_n)_{n}$ in $M$ defined  by (\ref{algoINC}), and satisfying  the following assumptions:
\begin{itemize} 

\item[$(i)$] For all $c>0$,
\[\sum_ne^{-c/\gamma_n} < \infty, \]

\item[$(ii)$] $(U_n)_n$ is uniformly bounded and
\[\mathbb{E} \left(U_{n+1} \mid \mathcal{F}_n \right) = 0,\]

\item[$(iii)$] $F$ is a standard set-valued map.
\end{itemize}

\hop Set $\tau_n := \sum_{i=1}^n \gamma_i$ and $m(t) := \sup \{ j \mid \tau_j \leq t \}$. We call $X$ the continuous time affine interpolated process induced by $(x_n)_n$ and $\overline{\gamma}$ the piecewise constant deterministic process induced by $(\gamma_n)_n$:
\[X (\tau_i + s) = x_i + s \frac{x_{i+1} - x_i}{\gamma_{i+1}}, \mbox{ for }  \; \, s \in [0,\gamma_{i+1}] \; { and } \; \;  \overline{\gamma}(\tau_i + s) := \gamma_{i+1} \mbox{ for } \; s \in [0, \gamma_{i+1}[,\]
and consider its limit set

\[\mathcal{L}(X) :=\bigcap_{t\ge 0}\overline{\{X(s) \ : \ s\ge t \}}.\]

\hop The following attainability condition is crucial to show that $X$ converges with positive probability to a given attractor.

\begin{definition} A point $p \in M$ is \emph{attainable} if, for any $t>0$ and any neighborhood $U$ of $p$,
\[\mathbb{P} \left(\exists s \geq t : \; \, X(s) \in U \right) >0.\]

\end{definition}

\hop We call $Att(X)$ the set of attainable points by $X$. The following statement is a special case of our main result, Theorem \ref{theorem}. 
\begin{theoreme}
Let $A \subset M$ be an attractor for $\Phi$ with basin of attraction $\mathcal{B}(A)$. If $Att(X) \cap \mathcal{B}(A) \neq \emptyset$ then
$$P(\mathcal{L}(X)  \subset A) > 0.$$
\end{theoreme}

\hop Heuristically, this means that if the set of attainable points of the process meets the basin of attraction of a given attractor $A$, then there is convergence with positive probability toward $A$.

\section{Convergence with positive probability}
\zcinq

\subsection{Set-valued dynamical systems relative to a differential inclusion.}
\ztrois

\hop We recall here some definitions and results due to Bena\"im et al (see \cite{BHS1}). 

\hop Let $F : \R^m \rightrightarrows \R^m$ be a standard set-valued map and $\Phi$ be the set-valued dynamical system associated to the differential inclusion
\begin{equation}\label{inclusiondiff2}
\frac{d \mathbf{x}}{dt} \in F(\mathbf{x})
\end{equation}

\hop The notion of {\itshape internally chain transitive set} (ICT set) was introduced by Bena\"im and Hirsch in \cite{BenHir96} to analyse certain perturbations of the flow relative to an ODE. This is an extension of the notion of chain recurrence due to Conley \cite{Con78}. The concept of ICT sets was extended to differential inclusions by Bena\"im et al. in \cite{BHS1}.

\hop We refer to this last reference for an accurate description of ICT sets. Here, we only need the following property (see \cite{BHS1}, Proposition 3.20): an internally chain transitive set $L$ is invariant, compact and the restricted set-valued dynamical system $\Phi \vert_L$ admits no proper attractor (i.e. no attractor distinct from $L$). The following result is proved in \cite{BHS1} (see Theorem 3.23) 

\begin{theoreme} \label{thmictattracteur} Let $L$ be an internally chain transitive set and $A$ be an attractor for $\Phi$ with basin of attraction $\mathcal{B}(A)$. Then
\[\mathcal{B}(A) \cap L \neq \emptyset \, \Rightarrow \; \, L \subset A.\]
\end{theoreme}
\zdeux

\hop The space $\mathcal{C}(\mathbb{R}_+, \mathbb{R}^m)$ of continuous paths, endowed with the metric 
\[\mathbf{D} (\mathbf{x},\mathbf{x'}) := \sum_{k=1}^{\infty} \frac{1}{2^k} \min \left( \sup_{u \in [0,k]} \|\mathbf{x}(u) - \mathbf{x'}(u) \|, 1 \right)\] 
is complete. A continuous map $X : \mathbb{R}_+ \rightarrow \mathbb{R}^m$ is an {\itshape asymptotic pseudotrajectory} (APT) of the set-valued dynamical system $(\Phi_t)_{t \geq 0}$ if, for any $T>0$,
\[\lim_{t \rightarrow + \infty} \inf_{\mathbf{z} \in S_{\Phi}} \left\|X(t+\cdot) - \mathbf{z}(\cdot) \right\|_{[0,T]} = 0,\]

\hop where $\|\cdot\|_{[0,T]}$ denotes the uniform norm on $[0,T]$. Heuristically this means that, for any $T>0$, the curve joining $X(t)$ to $X(t+T)$ shadows the trajectory of some solution  with arbitrary accuracy, provided $t$ is large enough.

\hop A fundamental property of asymptotic pseudotrajectories is that, if $X$ is a bounded APT, then its limit set $\mathcal{L}(X)$ is internally chain transitive (see \cite{BHS1}, Theorem 4.3). Consequently, by Theorem \ref{thmictattracteur}, we have

\begin{corollaire} \label{coraptattracteur} Let $X$ be an asymptotic pseudotrajectory of the set-valued dynamical system and $A$ an attractor for $\Phi$. If $\mathcal{L}(X)$ meets the basin of attraction of $A$, then it belongs to $A$.

\end{corollaire}

\hop Let $\delta$ be a positive real number. Then $F^{\delta}$ is the set-valued map defined by 
\begin{equation}\label{fdelta}
F^{\delta}(x) := \left\{y \mid \; \exists z \in B(x,\delta) \mbox{ such that } \; \, d(y,F(z)) < \delta \right\}.
\end{equation}

\hop 

\begin{definition} Let $\delta : \R_+^* \rightarrow \R_+^*$ and  $\overline{U} : \mathbb{R}_+ \rightarrow \R^m$. A continuous function $\mathbf{y} : \mathbb{R}_+ \rightarrow \R^m$ is a $(\delta(\cdot),\overline{U}(\cdot))$-perturbed solution of the differential inclusion (\ref{inclusiondiff2}) if 

\begitem

\item[$(i)$] $\mathbf{y}$ is absolutely continuous,

\item[$(ii)$] $\delta(t) \downarrow_{t \rightarrow + \infty} 0$ and,  for almost every $t > 0$, 
\[\frac{d\mathbf{y}(t)}{dt} - \overline{U}(t) \in F^{\delta(t)} (\mathbf{y}(t)),\]

\item[$(iii)$] $\overline{U}$  is locally integrable and such that, for any $T> 0$,
\begin{equation} \label{hypU}
\lim_{t \rightarrow + \infty} \sup_{0 \leq v \leq T} \left\|\int_t^{t+v} \overline{U}(s) ds \right\| = 0.
\end{equation}

\end{itemize}

\end{definition}

\hop We recall the following theorem due to Bena\"im et al (see \cite{BHS1} Theorem 4.2).

\begin{theoreme} \label{theoremPTA} Any bounded perturbed solution of the differential inclusion (\ref{inclusiondiff2}) is an asymptotic pseudotrajectory of the set-valued dynamical system $\Phi$.
\end{theoreme}
\zun

\subsection{A deterministic result}
\ztrois

For any application $X : \mathbb{R}_+ \rightarrow M$ and $T> 0$, we define the quantity
\[d_X(T) := \sup_{k \in \mathbb{N}} \inf_{\mathbf{z} \in S_{\Phi}}  \left\| \mathbf{z}(\cdot) - X(kT + \cdot) \right\|_{[0,T]}.\]

\hop The following characterization of the basin of attraction will be useful. 

\begin{lemme} Given an attractor $A$, the basin of attraction of $A$ is the union of all fundamental neighborhoods of $A$.
\end{lemme}

\hop {\bfseries Proof.} Any fundamental neighborhood is trivially included into $\mathcal{B}(A)$. Conversely, let $x \in \mathcal{B}(A)$ and $U_0$ be a given fundamental neighborhood of $A$: for any $\varepsilon > 0$, there exists $t^0_{\varepsilon} >0$ such that 
\[\Phi_t(U_0) \subset N^{\varepsilon}(A), \; \, \forall t \geq t^0_{\varepsilon}.\]

\hop Pick $\gamma>0$ such that $N^{2 \gamma}(A) \subset U_0$. There exists $T>0$ such that $\Phi_T(x) \subset N^{\gamma}(A)$. Otherwise, there would exist sequences $(\mathbf{z}_n)_n \subset S_{x}$ and $t_n \uparrow + \infty$ such that $\mathbf{z}_n(t_n) \in (N^{\gamma}(A))^c$, which contradicts the fact that $\omega_{\Phi}(x)$ is bounded and contained in $A$.
 
\hop Finally, by closedness of the set-valued map $x \mapsto \Phi_T(x)$\footnote{This is an easy consequence of Definition \ref{standarsetvaluedmap}.}, there exists $r>0$ such that $\Phi_T(B(x,r)) \subset N^{2\gamma}(A) \subset U_0$. The set $U := U_0 \cup B(x,r)$ is a fundamental neighborhood of $A$ since, for any $\varepsilon >0$, 
\[\Phi_t(B(x,r)) \subset N^{\varepsilon}(A), \; \, \forall t \geq t^0_{\varepsilon} + T, \]
and the proof is complete. $\; \; \blacksquare$
\zun

\hop Assume now that $X$ is an APT of the set-valued dynamical system $(\Phi_t)_{t \geq 0}$. The following lemma is the extension to differential inclusions of Lemma 6.8 in Bena\"im \cite{Ben99}.

\begin{lemme}\label{lemmeconvattracteur} Let $A \subset \R^m$ be an attractor for the set-valued dynamical system $\Phi$, with basin of attraction $\mathcal{B}(A)$. Then, for any compact set $K \subset \mathcal{B}(A)$, there exist positive real numbers $\alpha(K)$ and $T(K)$ such that
\[ \left( X(0) \in K \, \mbox{ and } \; \, d_X(T(K)) < \alpha(K) \right) \, \Rightarrow \mathcal{L}(X) \subset A.\]

\end{lemme}
\zun

\hop {\bfseries Proof.} Let $W$ be an open set with compact closure  such that 
\[A \cup K \subset W \subset \overline{W} \subset \mathcal{B}(A).\]

\hop There exists $\alpha > 0$ such that $N^{3 \alpha}(A) \subset W$ and $N^{\alpha}(K) \subset W$. Since $\overline{W}$ is included in the basin of attraction of $A$, there exists a fundamental neighborhood which contains $W$; hence we can find a positive number $T$ (which depends on $\alpha$ and $W$) such that
\[\Phi_{[T,+ \infty[}(W) \subset N^{\alpha}(A).\]

\hop Assume now that $X(0) \in K$ and $d_X(T) < \alpha$. There exists a solution $\mathbf{z^1}$ which shadows $X$ on $[0,T]$; in particular,
\[\mathbf{z^1}(0) \in N^{\alpha}(X(0)) \subset W \, \mbox{ and } X(T) \in N^{\alpha}(\mathbf{z^1}(T)). \]

\hop By definition of $T$, $\mathbf{z^1}(T) \in N^{\alpha}(A)$, which means that $X(T) \in N^{2 \alpha}(A) \subset W$.

\hop By a recursive argument, we show that the sequence $(X(kT))_{k \geq 1}$ belongs to the set $N^{2 \alpha}(A)$. Assume that $X(kT) \in N^{2 \alpha}(A)$. Then, there exists a solution $\mathbf{z^{k+1}}(\cdot)$ which is $\alpha$-close of $X(kT + \cdot)$ on $[0,T]$: in particular,
\begin{equation*}
\mathbf{z^{k+1}}(0) \in N^{\alpha}(X(kT)) \subset N^{3\alpha}(A) \subset W \, \mbox{ and } X(kT + T) \in N^{\alpha}(\mathbf{z^{k+1}}(T)) \subset N^{2\alpha}(A).
\end{equation*}

\hop Consequently, the limit set $\mathcal{L}(X)$ is contained in $W \subset \mathcal{B}(A)$. Hence, $\mathcal{L}(X) \subset A$ by Corollary \ref{coraptattracteur}. $\; \; \blacksquare$
\zdeux

\subsection{Stochastic processes}
\ztrois

\hop In the following,   $(X(t))_{t \geq 0}$ will be a  continuous time $\mathbb{R}^m$-valued stochastic process, adapted to some non decreasing sequence $(\mathcal{F}_t)_t$ of sub sigma algebras of $\mathcal{F}$.

\begin{hypothese}\label{H1} There exists a map $\omega : \mathbb{R}_+^3 \rightarrow \mathbb{R}_+$ such that, for any $\alpha > 0$ and $T > 0$, 
\begin{equation}\label{H12}
	\mathbb{P} \left(\sup_{s \geq t} \inf_{\mathbf{z} \in S_{\Phi}}\left\| \mathbf{z}(\cdot) - X(s + \cdot) \right		\|_{[0,T]} \geq \alpha \mid \mathcal{F}_t \right)  \leq \omega (t,\alpha,T) \downarrow_{t \rightarrow + 		\infty} 0 \; \mbox{ almost surely}.
\end{equation}

\end{hypothese}

\hop  Hypothesis \ref{H1} is a technical assumption, slightly stronger than supposing that $X$ is almost surely an APT. Recall that $Att(X)$ is the set of attainable points by $X$.

\begin{theoreme} \label{convposprob} Let $A$ be an attractor and $(X(t))_{t \geq 0}$ be an adapted process satisfying hypothesis \ref{H1}. Then, if $Att(X) \cap \mathcal{B}(A) \neq \emptyset$, we have
\[\mathbb{P} \left( \mathcal{L}(X) \subset A\right) > 0.\]

\end{theoreme}
\zun

\hop {\bfseries Proof.} We adapt the proof of Theorem 7.3 in Bena\"im \cite{Ben99}. Let $U$ be an open set included in $\mathcal{B}(A)$ and call $K = \overline{U}$. By Lemma \ref{lemmeconvattracteur} there exist $\alpha(K)$ and $T(K)$ such that
\[ (X(0) \in K \, \mbox{ and } \; \, d_X(T(K)) < \alpha(K)) \, \Rightarrow \mathcal{L}(X) \subset A.\]

\hop Let $t$ be a positive irrational number such that $\omega(t,\alpha,T) < 1$ and denote $t_n(k) = \frac{k}{2^n}$ for $n$ and $k$ in $\mathbb{N}$. We define the stopping time
\[\tau_n := \inf_{k \in \mathbb{N}} \left\{t_n(k) \mid X(t_n(k)) \in U, \; \, t_n(k) \geq t \right\}.\]

\hop On the intersection of the events  $\left\{\tau_n < \infty\right\}$ and  $\left\{\sup_{s \geq \tau_n} \inf_{\mathbf{z} \in S_{\Phi}} \|\mathbf{z}(\cdot) - X(s + \cdot)\|_{[0,T]} \leq \alpha \right\}$ the set $\mathcal{L}(X)$ is included in $A$. Consequently, we have
\begin{eqnarray*}
\mathbb{P} \left(\mathcal{L}(X) \subset A \right) &\geq& \sum_{k \geq [2^n t] + 1} \mathbb{E} \left( \mathbb{P} \left( \sup_{s \geq \tau_n} \inf_{\mathbf{z} \in S_{\Phi}} \|\mathbf{z}(\cdot) - X(s + \cdot)\|_{[0,T]} \leq \alpha \mid \mathcal{F}_{t_n(k)}\right) \mathbb{I}_{\tau_n = t_n(k)} \right) \\
&\geq& \sum_{k \geq [2^n t] + 1} \left(1 - \omega(t_n(k),\alpha,T)\right) \mathbb{P} \left(\tau_n = t_n(k) \right) \geq \left(1 - \omega(t,\alpha,T)\right) \mathbb{P} \left(\tau_n < + \infty \right),
\end{eqnarray*}

\hop since $\omega(t_n(k),\alpha,T) \leq \omega(t,\alpha,T), \, \; \forall k \geq [2^n t] +1$. On the other hand, the sequence of events 
$\{\tau_n < + \infty\}$ is increasing and we have
\[\lim_{n \rightarrow + \infty} \uparrow \{ \tau_n < + \infty\} = \cup_n \{ \tau_n < + \infty\} = \{\exists s \geq t \mid \; \, X(s) \in U\}.\]

\hop Now take an attainable point $p \in \mathcal{B}(A)$ and $U$ a neighborhood of $p$, such that $\overline{U} \subset \mathcal{B}(A)$. We have
\[\mathbb{P} \left( \mathcal{L}(X) \subset A\right) \geq (1 - \omega(t,\alpha,T)) \mathbb{P} \left(\exists s \geq t \mid \; \, X(s) \in U\right) > 0,\]

\hop and the proof is complete. $\; \; \blacksquare$
\zdeux

\hop Now, we consider a  compact set  $M \subset \R^m$ and a standard set-valued map $F : \R^m \rightrightarrows \R^m $. Let $T > 0$. Denote $\Phi^{T}(M) := \overline{\bigcup_{s\in [0,T]} \Phi_s(M)}$,  $||F|| = \sup_{x \in \Phi^T(M)} \sup_{y \in F(x)}||y||$ and let us define the compact set 
	$$K_C := K_{(||F||,C)} =  \left\{ \mathbf{y} \in Lip([0,T] , \R^m) \mid \; \, Lip(\mathbf{y}) \leq ||F|| + C + 1  \ ,  \   \mathbf{y}(0) \in M \right\},$$
where $C$ is a positive constant and $Lip([0,T],\R^m)$ is the set of Lipschitz functions from $[0,T]$ to $\R^m$.

\begin{Remarque}
$K_{C}$ is appropriate to our situation since it contains every solution curve, restricted to an interval of length $T$  and any $(\delta(\cdot),\overline{U}(\cdot))$-perturbed solution of the differential inclusion, with $\sup_{t \in [0,T]}\overline{U}(t)\le C$ and $\delta \le 1$.
\end{Remarque}	

\hop For  $\delta \in [0,1]$, let us define the set-valued application (with the convention $\Lambda^0 = \Lambda$):
\begin{equation}\label{lamdadelta}
\Lambda^{\delta} : K_C \rightrightarrows K_C, \; \, \mathbf{z} \mapsto \Lambda^{\delta}(\mathbf{z}),
\end{equation} 

\hop  where $\mathbf{y} \in \Lambda^{\delta}(\mathbf{z})$ if and only if there exists an integrable $h : [0,T] \rightarrow \R^m$ such that $h(u) \in F^{\delta}(\mathbf{z}(u)) \; \forall u \in [0,T]$ and 

\[\mathbf{y}(\tau) = \mathbf{z}(0) + \int_0^{\tau} h(u) du, \; \, \forall \tau \in [0,T].\]

\hop Remark  that $Fix(\Lambda) := \left\{\mathbf{z} \in K_C \mid \; \, \mathbf{z} \in \Lambda(\mathbf{z}) \right\}$ is the set of the restrictions  on $[0,T]$ of the solutions curves starting from $M$, which we denote $S_{[0,T]}$ from now on. Additionaly, we call  $d_{[0,T]}$ the distance associated to the uniform norm on $[0,T]$. The following lemma is an immediate consequence of Corollary \ref{distancefix}, proved in the last section.

\begin{lemme}\label{distFix} Let $C>0$ and $\alpha > 0$. There exist  $\varepsilon > 0$ and $\delta_0 > 0$ such that, for any $\delta < \delta_0$

$$d_{[0,T]}(\mathbf{z}, \Lambda^{\delta}(\mathbf{z}))< \epsilon \Rightarrow  d_{[0,T]} (\mathbf{z}, S_{[0,T]}) <\alpha.$$

\end{lemme}
\zdeux

\hop As a consequence, we obtain the following crucial result. 

\begin{proposition} \label{perturbedhyp} Assume that there exist a function $\delta : \mathbb{R}_+\rightarrow \mathbb{R}_+$ converging to zero and a uniformly bounded random process $(\overline{U}(t))_{t \geq 0}$ such that $(X(t))_{t \geq 0}$ is almost surely a bounded $(\delta,\overline{U})$-perturbed solution of the differential inclusion (\ref{inclusiondiff2}) and such that $X(0) \in M$. If $\overline{U}$ satisfies the following property

\begin{equation} \label{bruit}
\mathbb{P} \left(\sup_{s \geq t} \left\|\int_s^{s+\cdot} \overline{U}(u) du \right\|_{[0,T]} \geq \varepsilon \mid \mathcal{F}_t \right) \leq \omega(t,\varepsilon,T) \downarrow_{t \rightarrow + \infty} 0  
\end{equation}

\hop almost surely, then  hypothesis \ref{H1} holds and Theorem \ref{convposprob}  may be applied.

\end{proposition}
\zun

\hop {\bfseries Proof.} First, $X(\cdot)$ is almost surely an asymptotic pseudotrajectory of the set-valued dynamical system by Theorem \ref{theoremPTA}. For Lipschitz (classical) dynamical systems,  hypothesis \ref{H1} holds by an application of Gronwall lemma (see Bena\"im \cite{Ben99}, section 7). However, this does not adapt to our situation and this is the reason why we need Lemma  \ref{distFix}. By assumption, we have almost surely

\[\frac{d X(t)}{dt} - \overline{U}(t) \in F^{\delta(t)}(X(t)), \; \, \mbox{ for almost every } \, t>0. \]

\hop  Let $T>0$. For any $\tau \in [0,T]$, 

\[X(s+\tau) - \int_{s}^{s + \tau} \overline{U}(u) du \in  X(s) + \int_{0}^{\tau} F^{\delta(s)}(X(s+u)) du.\]

\hop Hence, $d_{[0,T]} \left( X(s + \cdot), \Lambda^{\delta(s)}(X(s+\cdot) \right) \leq  \left\|\int_s^{s + \cdot} \overline{U}(u) du  \right\|_{[0,T]}$ and

\begin{eqnarray*}
\mathbb{P} \left(\sup_{s \geq t} d_{[0,T]} \left(X(s+\cdot), \Lambda^{\delta(t)}(X(s+ \cdot))  \right)  \geq \varepsilon \mid \mathcal{F}_t \right) &\leq& \mathbb{P} \left(\sup_{s \geq t}   \left\|\int_s^{s+\cdot} \overline{U}(u) du \right\|_{[0,T]} \geq \varepsilon \mid \mathcal{F}_t \right) \\
&\leq& \omega(t,\varepsilon,T).
\end{eqnarray*}

\hop Now let $\alpha >0$. By Lemma \ref{distFix} there exists $\varepsilon > 0$ (which depends on $T$ and $\alpha$) such that, for $t$ large enough and $s \geq t$, 

\[d_{[0,T]} \left( X(s+\cdot), S_{[0,T]} \right) \geq \alpha \Rightarrow d_{[0,T]} \left(X(s + \cdot), \Lambda^{\delta(t)}(X(s + \cdot)) \right) \geq \varepsilon. \]

\hop Consequently, for these choices of $t$ and $\varepsilon$,

\begin{eqnarray*}
\mathbb{P} \left( \sup_{s \geq t} d_{[0,T]} \left( X(s + \cdot), S_{[0,T]} \right) \geq \alpha \mid \mathcal{F}_t \right) &\leq& \mathbb{P} \left(\sup_{s \geq t} d_{[0,T]} \left(X(s+\cdot), \Lambda^{\delta(t)}(X(s+ \cdot))  \right)  \geq \varepsilon \mid \mathcal{F}_t \right) \\
&\leq& \mathbb{P} \left( \sup_{s \geq t}  \left\|\int_s^{s+\cdot} \overline{U}(u) du \right\|_{[0,T]} \geq \varepsilon \mid \mathcal{F}_t\right)\\
&\leq& \omega(t,\varepsilon(\alpha,T),T),
\end{eqnarray*}

\hop and the proof is complete. $\; \; \blacksquare$
\zdeux

\hop Recall that a perturbed solution is an APT. The last proposition is a stochastic version of this result in the sense that the process $X(\cdot)$ is almost surely a perturbed solution, hence almost surely an APT. However, to yield the stronger property (\ref{H1}), we reinforce the assumption on $U$ and assume (\ref{bruit}) instead of (\ref{hypU}).

\subsection{Convergence of Stochastic approximation algorithms}
\ztrois

\hop We introduce here a class of stochastic approximation processes which generalize the Robbins-Monro algorithms. Under some assumptions on the step size and the noise, we prove that hypothesis \ref{H1} is verified and that the conclusion of  Theorem \ref{convposprob} holds.

\begin{definition}[Generalised stochastic approximation process]\label{robbinsmonro}

Let $(U_n)_n$ be an uniformly bounded $\mathbb{R}^m$-valued random process, $(\gamma_n)_n$ a deterministic positive real sequence and $(F_n)_{n}$  a sequence of set-valued maps on $\R^m$. We say that $(x_n)_n$ is a \emph{generalised stochastic approximation process relative to the standard set-valued map $F$} if the following assumptions are satisfied:

\begin{itemize} 

\item[$(i)$] we have the recursive formula

\[x_{n+1} -x_n - \gamma_{n+1} U_{n+1 }\in  \gamma_{n+1}  F_n(x_n), \]

\item[$(ii)$] the step size satisfies

\[ \sum_n \gamma_n = + \infty, \; \; \lim_n \gamma_n = 0,\]

\item[$(iii)$] for all $T > 0$, we have almost surely

\[\lim_{n \rightarrow + \infty} \sup \left\{ \left\| \sum_{i=n}^{k-1} \gamma_{i+1} U_{i+1} \right\| \mid \; \; k \, \mbox{ such that } \; \sum_{i=n}^{k-1} \gamma_i \leq T\right\} = 0,\]

\item[$(iv)$] for all $n\ge 0$, $x_n \in M$,

\item[$(v)$] for any $\delta > 0$, there exists $n_0 \in \mathbb{N}$
 such that
 
 \[\forall n \geq n_0, \; \, F_n(x_n) \subset F^{\delta}(x_n).\]

\end{itemize}

\end{definition}

\hop In the following we will call $X =(X(t))_t$ the continuous time affine interpolated process induced by a given generalised stochastic approximation process $(x_n)_n$ (see section 1.2).

\begin{proposition} \label{DSAperturbed} The process $X$ is almost surely a $(\delta,\overline{U})$-perturbed solution, for some deterministic function $\delta$, and $\overline{U}$ the piecewise constant continuous time process associated to $(U_n)_n$:
\[\overline{U}(t) := U_{n+1}, \; \; \forall t \in [\tau_n, \tau_{n+1}[.\]

\end{proposition}
\zun

\hop {\bfseries Proof.} By straightforward computations (see the proof of  proposition 1.3 in Bena\"im et al. \cite{BHS1}), it is not difficult to see that almost surely, $(X(t))_t$ is a perturbed solution associated to $\overline{U}$ and 
\[\delta (t) := \inf \left\{\delta > 0 \mid \; \; \tau_n \geq t \Rightarrow \, F_n(x_n) \subset F^{\delta}(x_n) \right\} + \overline{\gamma}(t) \left(\overline{U}(t) + c \left(1 + \sup_{x \in M} F(x) \right) \right),\]

\hop which obviously converges to $0$. $\; \; \blacksquare$

\begin{Remarque} Recall that $m(t) = sup \{j \mid \; \tau_j \leq t \}$. The condition (\ref{bruit}) is equivalent to 

\begin{equation} \label{bruit2}
\mathbb{P} \left(\sup_{m \geq n} \sup_{m < k \leq m(\tau_m +T)} \left\|\sum_{i=m}^{k-1} \gamma_{i+1} U_{i+1} \right\| \geq \varepsilon \mid \mathcal{F}_n \right) \leq \omega(n,\varepsilon,T) \downarrow_{n \rightarrow + \infty} 0 \; \mbox{ almost surely}.
\end{equation}

\hop and we use the notation  $\Delta(n,T) := \sup_{n < k \leq m(\tau_n +T)} \left\|\sum_{i=n}^{k-1} \gamma_{i+1} U_{i+1} \right\|$ in the sequel.

\end{Remarque}
\zdeux

\hop Our main result is now stated in full generality. 

\begin{theoreme}\label{theorem} Let $(x_n)_n$ be a stochastic approximation algorithm such that $(U_n)_n$ satisfies (\ref{bruit2}). Then, if $A$ is an attractor relative to $F$, we have

\[Att(X) \cap \mathcal{B}(A) \neq \emptyset \, \Rightarrow \; \, \mathbb{P} \left(\mathcal{L}(X) \subset A \right) > 0.\]

\end{theoreme}
\zun

\hop {\bfseries Proof.} By Proposition \ref{DSAperturbed}, the conditions requested to apply Proposition \ref{perturbedhyp} are satisfied. Hence, hypothesis \ref{H1} is checked and the result follows directly from Theorem \ref{convposprob}. $\; \; \blacksquare$ 
\zdeux

\hop In the particular case where $(U_n)_n$ is a martingale difference: $\mathbb{E} \left(U_{n+1} \mid \mathcal{F}_n \right) = 0$, (\ref{bruit2}) is satisfied under simple assumptions on the noise and step size.

\begin{proposition}\label{hypbruit} Let $(U^0_n)_n$ be a martingale difference noise (not necessarily bounded) and assume that one of the following assumptions is satisfied:

\begitem

\item[$1)$] There exists some $q \geq 2$ such that 
\[\sum \gamma_n^{1+ q/2} < + \infty \, \mbox{ and } \; \, \sup_n \mathbb{E}\left(\|U_n^0\|^q \right) < + \infty.\]

\item[$2)$] There exists a deterministic sequence $(M_n)_n$ such that $M_n^2 = o \left((\gamma_n \log n)^{-1}  \right)$ and, for any $n \in \mathbb{N}$, 

\[\forall \theta \in \mathbb{R}^m, \; \, \mathbb{E} \left(\exp \left( \left< \theta, U^0_{n+1} \right> \right) \mid \mathcal{F}_n \right) \leq \exp \left( \frac{M_n^2}{2} \|\theta\|^2 \right),\]

\end{itemize}

\hop then (\ref{bruit2}) is checked.

\end{proposition}

\hop {\bfseries Proof.} For the first point, we refer the reader to Bena\"im \cite{Ben99}. Now for the second, let $\theta \in \mathbb{R}^m$ and consider the process $(Z_n(\theta))_n$ defined by

\[Z_n(\theta) := \exp \left( \sum_{i=1}^n \left<\theta, \gamma_i U^0_i \right> - \frac{\|\theta\|^2}{2} \sum_{i=1}^n \gamma_i^2 M_i^2\right). \]

\hop $(Z_n(\theta))_n$ is a supermartingale by assumption. Hence, if we denote $S_n := \sum_{i=n}^{k-1} \gamma_{i+1}^2 M_{i+1}^2$ and $m_n := m(\tau_n + T)$, for any $\beta > 0$, 

\begin{eqnarray*}
\mathbb{P} \left(\sup_{n < k \leq m_n} \left<\theta, \sum_{i=n}^{k-1} \gamma_{i+1} U_{i+1} \right> \geq \beta \mid \mathcal{F}_n  \right) &=& \mathbb{P} \left( \sup_{n < k \leq m_n} \frac{Z_k(\theta)}{Z_n(\theta)} \exp \left( \frac{\|\theta\|^2}{2}\sum_{i=1}^k \gamma_i^2 M_i^2\right) \geq e^{\beta} \mid \mathcal{F}_n   \right) \\
&\leq&  \mathbb{P} \left( \sup_{n < k \leq m_n} Z_k(\theta) \geq Z_n(\theta) \exp \left(\beta - \frac{\|\theta\|^2}{2} S_n \right) \mid \mathcal{F}_n   \right) \\
&\leq& \exp \left(\frac{\|\theta\|^2}{2} S_n - \beta \right).
\end{eqnarray*}

\hop Let $e \in \left\{e_1,...,e_m,-e_1,...,-e_m \right\}$. We have

\begin{eqnarray*}
\mathbb{P}\left(\sup_{n < k \leq m_n} \left<e, \sum_{i=n}^{k-1} \gamma_{i+1} U_{i+1} \right> \geq \varepsilon \mid \mathcal{F}_n  \right) &=& \mathbb{P}\left(\sup_{n < k \leq m_n} \left<\frac{\delta e}{S_n}, \sum_{i=n}^{k-1} \gamma_{i+1} U_{i+1} \right> \geq \frac{\varepsilon^2}{S_n} \mid \mathcal{F}_n  \right) \\
&\leq& \exp \left( \frac{-\varepsilon^2}{2 S_n}\right).
\end{eqnarray*}

\hop 

\[\mathbb{P}\left(\sup_{n < k \leq m_n} \left\|\sum_{i=n}^{k-1} \gamma_{i+1} U_{i+1} \right\| \geq \varepsilon \mid \mathcal{F}_n  \right) \leq 2m \exp \left( \frac{-\varepsilon^2}{2 S_n}\right).\]

\hop Let us introduce $\varepsilon_n := \gamma_n M^2_n \log n$.  Then, since $\sum_{i=n}^{k-1} \gamma_{i+1} \leq T$, we have

\[\mathbb{P}\left(\sup_{n < k \leq m_n} \left\|\sum_{i=n}^{k-1} \gamma_{i+1} U_{i+1} \right\| \geq \varepsilon \mid \mathcal{F}_n  \right) \leq 2m \exp \left( \frac{-\varepsilon^2 \log n}{2T \sup_{k \geq n} \varepsilon_k}\right).\]

\hop Since $\sup_{k \geq n} \varepsilon_k \rightarrow 0$, the application $(n,\varepsilon,T) \mapsto \omega(n,\varepsilon,T) := 2d \sum_{m \geq n} \exp \left(\frac{-\varepsilon^2 \log n}{2T \sup_{k \geq n} \varepsilon_k} \right)$ converges to $0$ as $n$ tends to infinity and the proof is complete. $\; \; \blacksquare$ 
\zdeux

\hop Simple examples satisfying the assumptions of Proposition \ref{hypbruit} are $\gamma_n = 1/n$ and $(U^0_n)_n$ a martingale difference with uniformly bounded moment of order $2$ (for $1)$) or $\gamma_n = (1/log(n))^2$ and $(U^0_n)_n$ a uniformly bounded martingale difference (for $2)$).
\zdeux

\hop The following is a useful consequence of this statement.

\begin{corollaire}
Assume that $U_n$ can be written $U_n = U_{n}^{0}+U_{n}^{1}$, where
\begin{itemize}
\item $ (U_{n}^{0})_n$ a martingale difference noise satisfying one of the assumptions in the Proposition \ref{hypbruit},
\item (\ref{bruit2}) is satisfied for $(U_{n}^{1})_n$.
\end{itemize}
Then (\ref{bruit2}) is satisfied for $(U_{n})_n$
\end{corollaire}

\hop {\bfseries Proof.} The sum of two random sequences satisfying (\ref{bruit2})  enjoys the same property. $\; \; \blacksquare$

\section{Application to the Markovian fictitious play learning model}
\zcinq

\hop We discuss here a Markovian strategy in a two-person game and study the induced dynamics. The model is studied by Bena\"im and Raimond in \cite{BenRai08} and was inspired by a so-called {\itshape pairwise comparison dynamics} introduced in Bena\"im et al. \cite{BHS2}.  

\subsection{The model}
\ztrois

\hop The motivation is the following. We assume, in the initial model, that the \emph{information situation} is the same as in the smooth fictitious play developped by Fudenberg and Levine (see \cite{FudKre93} and  \cite{FudLev98}) where the considered player uses a best response strategy against the empirical moves of his opponent, with respect to a smooth perturbation of the payoff function. A player adopting a smooth fictitious play strategy needs to be informed of his payoff function as well as the moves of his opponents up to this stage.

\hop For some reason (for instance if his set of actions is too large, if he has computational limitations or, more simply, if he is not allowed to play every action at each stage), we consider here that the set of moves he can play at some instant is a subset of his action set, which depends on the last action taken. 

\hop More formally, we consider a two players game in normal form. Let $I$ and $L$ be the (finite) sets of moves of respectively player 1 and player 2. These sets are of the form
\[I = \{1,...,m^1\}, \; \; L = \{1,...,m^2\};\]

\hop The maps $(U^1,U^2) : I \times L \rightarrow \mathbb{R} \times \mathbb{R}$ denote the payoff (or utility) functions of players. The sets of mixed strategies available to players are denoted $\mathcal{X} = \Delta(I)$ and $\mathcal{Y} = \Delta(L)$, where
\[\Delta(I) := \left\{ x=(x_1,...,x_{m^1}) \in \mathbb{R}_+^{m^1} \, \mid \, \; \sum_{i=1,..,m^1} x_i =1 \right\},\]

\hop and analogously  for $\Delta(L)$. The product $\mathcal{X} \times \mathcal{Y}$ is denoted $\Delta$. We will use the classical abuse of language for $y \in \mathcal{Y}$:
\[U^1(i,y) = \sum_{l \in L} U^1(i,l) y_l.\]

\hop For $x \in \mathcal{X}, \; y \in \mathcal{Y}$, we call $br^1(y) := Argmax_{x \in \mathcal{X}} U^1(x,y)$ and $br^2(x) = Argmax_{y \in \mathcal{Y}} U^2(x,y)$. We define the set-valued map $F :\mathcal{X} \times \mathcal{Y} \rightrightarrows \mathcal{X} \times \mathcal{Y}$  by 
$$F(x,y) = \{(\alpha,\beta)  \, \mid \; \, \alpha \in br^1(y) -x, \; \beta \in br^2(x) -y\}.$$
We assume that a given game is played repeatedly and call $X_n$ (resp. $Y_n$) the move of player 1 (resp. player 2) at stage $n$. The empirical distribution of moves up to stage $n$ is denoted $\overline{x}_n$ (resp. $\overline{y}_n$). We 
\zdeux

\hop Let $M_0^1$ be an irreducible matrix, reversible with respect to its invariant probability distribution $\pi^1_0$, which means that 
\[(\pi^1_0)_i M_0^1(i,j) = (\pi^1_0)_j M^1_0(j,i).\]
The matrix $M^1_0$ represents the possibility or not to play an action depending on the last move: player 1 will be able to play action $j$ after having played $i$ if and only if $M^1_0(i,j) >0$. For $n \in \mathbb{N}$ and $y \in \mathcal{Y}$, let us define the Markov matrix
\[M^1_n(i,j;y) = \left\{ 
            \begin{array}{ll}
               M^1_0(i,j) \exp \left( -\beta^1_n \left( U^1 \left(i,y \right) - U^1 \left(j, y\right)\right)^+\right) \; & \mbox{if} \, i \neq j, \\
               1 - \sum_{k \neq i} M^1_n(i,k;y) & \mbox{if} \, i=j,
            \end{array}
            \right.
\]
where $(\beta^1_n)_n$ is some positive deterministic sequence.

\begin{definition} A Markovian fictitious play (MFP) strategy for player $1$, associated with $(\beta^1_n)_n$ and $(M^1_0,\pi^1_0)$ is a strategy $\sigma$ such that, for any $n \in \mathbb{N}$, 
\[\mathbb{P}_{\sigma} \left(X_{n+1}=j \mid \mathcal{F}_n \right) = M^1_n(X_n,j;\overline{y}_n). \]
\end{definition}

\hop From now, we assume that both players use a Markovian fictitious play strategy, associated to $M_0^p$ and $(\beta_n^p)_n$ ($p=1,2$). Let us introduce the random sequences
\[V_n := \left(\delta_{X_n}, \delta_{Y_n} \right) \, \mbox{ and } \; \, v_n := \frac{1}{n} \sum_{i=1}^n V_i = \left(\overline{x}_n, \overline{y}_n \right).\]

\hop We call $Att(v)$ the attainability set of the discrete process $(v_n)_n$. Recall that $p \in Att(v)$ if and only if, for any neighborhood $N$ of $p$ and any $n_0 \in \mathbb{N}$, 

\[\mathbb{P} \left( \exists  n \geq n_0 : \; v_n \in N\right) > 0.\]

\begin{proposition} Assume that the matrices $M_0^p$ ($p=1,2$) have positive diagonal entries. Then $Att(v)$ is equal to the whole state space $\Delta$.
\end{proposition}
\zun

\hop {\bfseries Proof.} By irreducibility, from any instant $n$, given any player $i$, any move $a^i_j$ and any positive integer $p$, player $i$ will play action $a^i_j$ $p$ times in a row with positive probability. $\; \; \blacksquare$

\begin{theoreme} \label{ConvposMFP} Assume that the matrices $M_0^p$ ($p=1,2$) have positive diagonal entries. There exist positive values $\tilde{A}^p, \; p=1,2$ (which depend only on the payoff functions and $M_0^p$, $p=1,2$) such that, if agent $p$ plays accordingly to a MFP strategy with $\beta^p_n = A^p \log n$ and $A^p < \tilde{A}^p$, then 
\[\mathbb{P} \left(\mathcal{L}((v_n)_n) \subset A \right) > 0,\]
for any attractor  $A$ for the set-valued dynamical system induced by $F$.
\end{theoreme}

\hop In particular, a strict Nash equilibrium is always an attractor for the best response dynamics. Hence

\begin{corollaire} Let $\hat{v} = (\hat{x}, \hat{y})$ be a strict Nash equilibrium. Then, under the assumptions of Theorem \ref{ConvposMFP}, 
\[\mathbb{P} \left( v_n \rightarrow \hat{v} \right) > 0.\]
\end{corollaire}

\hop Another consequence of Theorem \ref{ConvposMFP} is the following. Assume that $U^1 = U^2 = U$ (we will call such a game a \emph{potential game}) and call $\Lambda$ the set of local maximizers of $U$:
\[\Lambda = \{u \in \Delta : \exists  V_u \in \mathcal{N}_u : \forall v \in V_u, U(v) \le U(u) \},\]
where $\mathcal{N}_u$ is the set of open neighborhoods of $u$. The set of Nash equilibria is denoted by $NE$.

\begin{corollaire} \label{potential} Assume that $L$ is a closed (in $\Delta$) connected component of $\Lambda$, which is isolated in the sense that there exists an open neighborhood $W$ of $L$ such that $W \cap \Lambda = L$. Then $L$ is an attractor for the best-response dynamics  and therefore, under the assumptions of Theorem \ref{ConvposMFP},
\[\mathbb{P}(\mathcal{L}((v_n)_n) \subset L) > 0.\]
\end{corollaire}
\hop {\bfseries Proof.} The fact that $L$ is an attractor for the best response dynamics is proved in greater generality in section \ref{proofPotential}, Proposition \ref{propatt}. $\; \; \blacksquare$

\begin{Remarque}
Notice that the closedness $\Lambda$ is essential. Consider a simple example where both players have two actions and the common payoff matrix is 
$$ A= \begin{pmatrix}
1 & 1 \\
0 & 2
\end{pmatrix}$$
A direct computation shows that 
\[NE = \{ \big((1,0), (1-t,t)\big) : \  t \in [0, 1/2]\} \cup \{\big((0,1),(0,1)\big)\}.\]
\hop Let us call the first set $L_1$ and the second $L_2$. We easily see that $L_2$ is an attractor for the best response dynamics. On the other hand, $L_1$ is a closed connected component of $NE$ but not of $\Lambda$  since the Nash equilibrium $(x,y) = \big((1,0),(\frac{1}{2},\frac{1}{2})\big)$ is not a local maximizer. Thus the connected component $L_1 \backslash \{(x,y)\}$ of $\Lambda$ is not isolated in $NE$ and Corollary \ref{potential} can not be applied as it will be made clear in its proof.
\end{Remarque}

\subsection{Proof of Theorem \ref{ConvposMFP}} \label{proofMFP}
\ztrois

\hop Notice that the Markov matrix $M^1_n(\cdot,\cdot;y)$ defined in the  previous section is reversible with respect to its invariant distribution $\pi^1_n[y]$:
\[ (\pi^1_n[y])_i \propto  (\pi^1_0)_i \exp \left(\beta_n U^1(i,y) \right).\]
Also, considering an irreducible Markovian matrix $M$ and its invariant probability measure $\pi$, one can define the pseudo inverse  $Q$ of $M$,  characterized by
\[Q(I-M) = (I-M)Q = I - \Pi, \; \; Q \, \mathbf{1} = 0, \]

\hop where $\Pi$ is the matrix defined by $\Pi(i,j) = \pi(j)$. Let us call $\pi^1_n$ and $Q^1_n$ (respectively $\pi_n^2$ and $Q_n^2$) the invariant probability and the pseudo inverse of the matrix $M^1_n := M^1_n(\cdot,\cdot;\overline{y}_n)$ (resp. $M_n^2 := M_n^2(\cdot,\cdot; \overline{x}_n)$). We now define the {\em energy barrier} of  $M^1_0$ with respect to the payoff function $U^1$. Let $\Gamma_{i,j}$ be the set of admissible paths from $i$ to $j$ in the graph associated to $M^1_0$: $\gamma = (i=i_0,i_1,..,i_n=j)$ is admissible if $M^1_0(i_k,i_{k+1}) > 0, \; k=0,..,n-1$. Then, denoting for $y \in \mathcal{Y}$, 
\[Elev(i,j;y) := \min \left\{ \max \{-U^1(k,y) \mid \, k \in \gamma\}, \; \gamma \in \Gamma_{i,j} \right\}.\]
We call
\[U^{1,\#}(y) :=  \max \left\{Elev(i,j;y) + U^1(i,y) + U^1(j,y) - \max U^1(\cdot,y) \right\}, \; \; U^{1,\#} := \max_{y \in \mathcal{Y}} U^{1,\#}(y).\]
Obviously, the quantity $U^{2,\#}$ is defined analogously. 

\hop For $v=(x,y) \in \Delta$, we call $\theta_n(v)$ the random variable
\[\theta_n(v) := \left(\pi^1_n[y], \pi^2_n[x]  \right).\]
The stochastic process $(v_n)_n$ satisfies the recursive formula
\[v_{n+1} - v_n = \frac{1}{n+1} \left(-v_n + V_{n+1} \right) = \frac{1}{n+1} \left(-v_n + \theta_n(v_n) + U_{n+1} \right),\]
with
\[U_{n+1} = V_{n+1} - \theta_n = \left(\delta_{X_{n+1}} - \pi^1_n[\overline{y}_n], \delta_{Y_{n+1}} - \pi^2_n[\overline{x}_n]  \right).\]
\zun

\hop  Bena\"im and Raimond (see Theorem 4.15 in \cite{BenRai08}) proved that, if $\beta_n^p = A^p \log n$ with $A^p < \tilde{A}^p := 1/2 U^{p,\#}$, then  $(v_n)_n$ is a generalised stochastic approximation process, taking values in $\Delta$, with step size $\gamma_n = 1/n$,  relatively to the maps $F_n(v) = -v + \theta_n[v]$ and  $F$. Note that the corresponding differential inclusion is the best response dynamics:
\[(\dot{x},\dot{y}) \in (br^1(y), br^2(x)) - (x,y)\]

\hop However, this is not sufficient to prove Theorem \ref{ConvposMFP} and we need to state the stronger property (\ref{bruit2}) on $(U_n)_n$. 

\hop The following proposition can be easily derived from the proof of Proposition 4.4 in Bena\"im and Raimond \cite{BenRai08}. 

\begin{proposition} \label{beta} Assume that $\beta^1_n = A^1\log n$   for some $0 < A^1 < \tilde{A}^1$. Then there exists a positive deterministic sequence $(u_n)_n \rightarrow 0$ such that

\begitem

\item[$a)$] $\frac{|Q^1_n|^2 \log n}{n} \leq u_n$,

\item[$b)$] $ |\Pi^1_{n+1} - \Pi^1_n| \leq u_n$,

\item[$c)$] $ |Q^1_{n+1} - Q^1_n| \leq u_n$.

\end{itemize}

\end{proposition}

\begin{lemme} Assume that the sequences $(\beta_n^p)_n$ ($p=1,2$) satisfy the assumption of Proposition \ref{beta}.  Then (\ref{bruit2}) is satisfied for  $(1/n)_n$ and $(U_n)_n$. 
\end{lemme}
\zun

\hop {\bfseries Proof.}   We call $\zeta_{n+1}$ the term $\delta_{X_{n+1}} - \pi^1_n[\overline{y}_n]$. We only need to prove that property (\ref{bruit2}) holds for $1/n$ and $(\zeta_n)_n$. We therefore denote $\Delta (n,T) := \sup_{n < k < m(\tau_n + T)} \left\| \sum_{i=n}^{k-1} \frac{1}{i+1} \zeta_{i+1}\right\|$.

\hop First of all, $\zeta_{n+1}$ can be written
\[\zeta_{n+1} = \delta_{X_{n+1}} \left(Id - \Pi_n \right) = \delta_{X_{n+1}} \left(Q_n - M_n Q_n \right).\]

\hop There is then a natural decomposition:
\[\zeta_{n+1} = (\delta_{X_{n+1}} Q_n - \delta_{X_n} M_n Q_n) + (\delta_{X_n} M_n Q_n - \delta_{X_{n+1}} M_n Q_n).\]

\hop The first term is a martingale difference, bounded by $|Q_n|$ (up to a constant). Hence it satisfies the assumption $2)$ of Proposition \ref{hypbruit}, with $M_n = \sqrt{\frac{n u_n}{\log n}}$.

\hop Now, for the second term, we have
\begin{eqnarray*}
\sum_{i=n}^{k-1} \frac{1}{i+1} \left(\delta_{X_i} M_i Q_i - \delta_{X_{i+1}} M_i Q_i \right) &\leq& \sum_{i=n}^{k-1} \frac{1}{i+1}\left(\delta_{X_{i+1}} M_{i+1} Q_{i+1} - \delta_{X_{i+1}} M_i Q_i \right) \\
&+& \sum_{i=n}^{k-1} \left(\frac{1}{i} \delta_{X_i} M_i Q_i - \frac{1}{i+1} \delta_{X_{i+1}} M_{i+1} Q_{i+1}\right) + T \sup_{n \leq i \leq k-1} \frac{|Q_i|}{i},
\end{eqnarray*}
since $\left\|\sum_{i=n}^k \frac{1}{i(i+1)} \delta_{X_i} M_i Q_i \right\| \leq \sup \left\{|Q_i|/i \mid n \leq i \leq k-1\right\}  \sum_{i=n,..,k-1} 1/i$.

\hop The first term on the right side can be written
\[\sum_{i=n}^{k-1} \frac{1}{i+1} \delta_{X_{i+1}} \left( Q_{i+1} - Q_i + \Pi_{i+1} - \Pi_i \right),\]
\hop and is bounded by the quantity $T \sup \left\{|Q_{i+1} - Q_i| + |\Pi_{i+1} - \Pi_i| \, \mid \; \, i=n,..,k-1 \right\}$. The telescopic term is bounded by $2 \sup_{n \leq i \leq k-1} \frac{|Q_i|}{i}$. Consequently, 
\[\left\| \sum_{i=n}^{k-1}  \frac{1}{i+1} \left(\delta_{X_i} M_i Q_i - \delta_{X_{i+1}} M_i Q_i \right)\right\| \leq 2(T+1) \sup_{i \geq n} \left\{\frac{|Q_i|}{i} + |Q_{i+1} - Q_i| + |\Pi_{i+1} - \Pi_i|\right\}. \]

\hop By Proposition \ref{beta}, the term on the right is decreasing to zero and
\[\mathbb{P} \left( \sup_{m \geq n} \Delta(n,T) \geq  \varepsilon \mid \mathcal{F}_n \right) \leq \omega(n,\varepsilon,T) \downarrow_n 0.\]

\hop This concludes the proof.  $\; \; \blacksquare$

\subsection{Potential Games and Proof of Corollary \ref{potential}}\label{proofPotential}
\ztrois

\hop In this section, we prove that any closed isolated connected component of $\Lambda$, the set of local maximizers, is an attractor for the best response dynamics in any given finite $N$-players potential game. Let $(m^i)_{i=1}^N$ be $N$ natural numbers and assume that the set of pure action for player $i$ is $\{1,..,m^i\}$. We call $\Delta^{m^i}$ the $m^i-1$-dimensional simplex corresponding to its mixed strategies space and $U$ the common $n$-linear payoff function on $\Delta = \times_{i=1}^N \Delta^{m^i}$ (see section 3.1). We use the notation $(x^i, x^{-i})$  for $x = (x^1,...,x^n) \in \Delta$. 
\ztrois

\begin{Remarque}
The set $\Delta$ can be written as a finite union of the relative interiors $(F_k)_{k=1..K}$ of its faces (this is also true for more general convex sets, see R. Tyrrell Rockafellar \cite{RockTyr}, Theorem 18.2). Additionaly, any $F_k$ is equal to a product $F_k^1 \times \cdots \times F_k^N$, where $F_k^i$ is the relative interior of one of $\Delta^{m^i}$'s faces. The restriction of $U$ to $F_k$ can therefore be seen as a smooth function defined on an open set of $\mathbb{R}^{n_k}$, for some natural number $n_k$.
\end{Remarque}

\hop The best response map $BR(x) = (BR^1(x), \dots, BR^n(x))$, where $BR^i(x) = Argmax_{z\in \Delta^{m^i}}U(z,x^{-i})$, has nonempty compact convex values and is upper semicontinuous. Therefore  consider the best response dynamics
$$\stackrel{.}{x} \  \in  Br(x) -x.$$
Bena\"im, Hofbauer and Sorin, in \cite{BHS1}, proved that $U$ is a Lyapunov function for $NE$. Namely, $U$ verifies the following two properties
\begin{itemize}
 \item[(i)] $U(x) < U(y)$ for all $x \in \Delta \backslash NE$, $y \in \Phi_t(x)$, $t>0$,
 \item[(ii)] $U(x) \le U(y)$ for all $x \in NE$, $y \in \Phi_t(x)$, $t\ge0$. 
\end{itemize}

\begin{lemme}\label{lemmeNE}
Let $\Lambda$ be the set of local maximizers of $U$. Then we have
$$ \Lambda \subset NE.$$
\end{lemme}
\hop {\bfseries Proof.}
Let $x$ be a local maximizer of $U$. Since $U(\cdot, x^{-i})$ is linear for all $i$, we easily derive that $x^i$ is a global maximizer for $U(\cdot, x^{-i})$. $\; \; \blacksquare$

\begin{proposition}\label{propNbCompCon}
For any finite game $\Gamma$, the set of Nash equilibria consists of finitely many connected components.
\end{proposition}
\hop For a proof of the above result see Kohlberg and Mertens \cite{KohlMert86}.

\begin{lemme}\label{lemmecritic}
Let $k \in \{1,..,K\}$. Then
$$ NE \cap F_k \subset \Sigma_{U_{|F_k}},$$
where $\Sigma_{U_{|F_k}}$ is the set of critical points of $U_{|F_k}$, the restriction of $U$ to $F_k$.
\end{lemme}

\hop {\bfseries Proof.}
Let $\tilde{x}$ be an element of $NE \cap F_k$. By definition of Nash equilibria, we know that, for all $i \in \{1, \dots, N\}$, the point $\tilde{x}^i$ is a global maximizer in $F^i_k$ for the function $U_{|F_k}(\cdot,\tilde{x}^{-i})$. In particular, $\tilde{x}$ is a critical point for $U_{|F_k}$. $\; \; \blacksquare$ 

\begin{lemme}\label{lemmeConstant}
The potential function $U$ is constant on every connected component $C$ of $NE$.
\end{lemme}
\hop {\bfseries Proof.}
Let $C$ be a connected component of $NE$ in $\Delta$. First, we write $C = \cup_{k =1..K} C_k$, where, for any $k$, $C_k = C \cap F_k$. Pick $k \in \{1,..,K\}$.  Since $U_{|F_k}$ is smooth, the image of its critical set has null Lebesgue measure by Sard's Lemma. Hence, by Lemma \ref{lemmecritic}, this is also true for the image of $C_k$. In particular, $U(C)$ has null Lebesgue measure, which means that it contains no interval. The result follows by connectedness of $C$ and continuity of $U$.
\hop  $\; \; \blacksquare$

\begin{Remarque} This result still holds for the more general class of potential games studied in \cite{BHS1}, where the strategy spaces are convex compact subsets of euclidean spaces, with countably many faces and the common payoff function is just assumed to be $N$-concave and smooth. However, Proposition \ref{propNbCompCon} is not, a priori, true anymore.
\end{Remarque}

\begin{lemme}\label{compconLNE}
Assume that $L$ is a closed (in $\Delta$) isolated connected component of $\Lambda$. Then it is also an isolated connected component of $NE$.
\end{lemme}
\hop {\bfseries Proof.}
By Lemma \ref{lemmeConstant}, there exists a constant $c \in \mathbb{R}$ such that $U(L) = \{c\}$. Assume that there exists a sequence $(x_n) \subset NE \backslash L$ and $x \in L$ such that $(x_n)$ converges to $x$. By  Proposition \ref{propNbCompCon}, there exists a  connected component of Nash equilibria $C$ which contains $L$ and  such that, for $n$ large enough, $x_n \in C$. In particular, $U(x_n) =c$ by Lemma \ref{lemmeConstant}. Moreover, since $x$ is a local maximizer, there is a neighborhood $V$ of $x$ such that for all $y \in V$, $U(y) \le c$. For $n$ large enough, $x_n \in V$ and it is also a local maximizer which is a contradiction. Therefore $L$ is isolated in the set of Nash equilibria and $L=C$. $\; \; \blacksquare$
\zun

\hop  The following proposition is a simplified version of Proposition 3.25 in \cite{BHS1}.
\begin{proposition}\label{proplyapunov}
Let $M$ be an Euclidean space, $K \in M$ be a compact set,  $V \subset  M$ be a bounded open neighborhood of $K$  and $U : \overline{V} \longrightarrow \R_-$. Let the following hold :
\begin{description}
\item[$a)$] For all $t\ge0$, $\Phi_t(V) \subset V$ (i.e., $V$ is strongly positively invariant);
\item[$b)$] $U^{-1}(0) = K$;
\item[$c)$] $U$ is continuous and for all $u \in V \setminus K$, $t>0$ and $v \in \Phi_t(u)$, $U(v) > U(u)$.
\end{description}
Then there exists an attractor contained in $K$ whose basin contains $V$, and with $U^{-1}(]-r,0])$ as a fundamental neighborhood for $r>0$ small enough.
\end{proposition}

\begin{proposition}\label{propatt}
If $L$ is a closed (in $\Delta$) isolated connected component of  $\Lambda$, then it is an attractor for the best response dynamics.
\end{proposition}

\hop {\bfseries Proof.}
By Lemma \ref{lemmeConstant}, there exists a number $c \in \R$ such that $L \subset U^{-1}(c)$. Without loss of generality we assume that $c=0$. 

\hop In order to use Proposition \ref{proplyapunov}, we construct an appropriate open neighborhood $V_r$ of $L$. First,  by Lemma \ref{compconLNE}, there exists an open neighborhood $W$ of $L$ such that
\begin{equation} \label{isolation}
W \cap NE = L.
\end{equation} 
Moreover there exists another open neighborhood $V$ of $L$ such that $\overline{V} \subset W$ and, for all $u \in \overline{V} \backslash L$, $U(u) < 0$. Indeed, assume by contradiction that there exists a sequence $(u_n)_n \subset W \backslash L$ which converges to $u \in L$ and such  that $U(u_n) \ge 0$ for all $n \ge 0$. There exists a neighborhood $V_u$ of $u$ such that for all $v \in V_u$, $U(v) \le 0$. For $n$ large enough, $V_u$ is also a neighborhood of $u_n$, hence $U(u_n) = 0$ and $u_n$ belongs to $\Lambda$, which contradicts (\ref{isolation}).

\hop By continuity of $U$, there exists a real number $r$ such that $U^{-1}(]-r,+\infty[) \cap \overline{V} \subset V$. Pick $V_r = U^{-1}(]-r, + \infty[ \cap \overline{V}$ (notice that $V_r$ is an open set, included in $V$) and consider the function $\tilde{U} : \overline{V}_r \longrightarrow \R$ which is the restriction  of the function $U$ to $\overline{V}_r$. 


\hop By construction of $V_r$, we clearly have $\tilde{U}^{-1}(0) = L$ and point $b)$ is therefore checked. Now let $t > 0$, $u \in V_r \backslash L$ and $v \in \Phi_t(u)$. By Property $(i)$ of a Lyapunov function (see p.16), $U(v) > U(u) > -r$ and point $c)$ is satisfied. Now assume that $v \notin V_r$. Then, since $v \in U^{-1}(]-r,+ \infty[)$,  $v \notin \overline{V}$ and there exists some $0< t'<t$ and $w \in \Phi_{t'}(u) \cap (\overline{V} \setminus V_r)$, which implies that $U(u) < U(w) \leq -r < U(u)$, a contradiction. Consequently $v \in V_r$. This proves that points $a)$ is checked and Proposition \ref{proplyapunov} applies: there exists an attractor included in $L$ and whose basin of attraction contains $V_r$. The proposition is proved since $L$ clearly admits no proper attractor. $\; \; \blacksquare$

\section{Proof of Lemma \ref{distFix}}
\zcinq

\hop In the following, $L^1([0,T])$ is the set of all Lebesgue-integrable functions from $[0,T]$ to $\mathbb{R}^m$. Let $H : [0,T] \rightrightarrows \mathbb{R}^m$ be a set-valued map, such that, for any $u \in [0,T]$, $H(u)$ is a nonempty subset of $\mathbb{R}^m$.

\begin{definition} We call $\mathcal{S}(H)$ the set of integrable selections from $[0,T]$ to $\mathbb{R}^m$:

\[\mathcal{S}(H) := \left\{h \in L^1([0,T]) \mbox{ such that } \; \, \forall u \in [0,T], \; \,  h(u) \in H(u) \right\}. \]

\hop With such a definition, we introduce the set-valued integral of $H$ on $[0,T]$:

\[\int_{[0,T]} H(u) du := \left\{ \int_{[0,T]} h(u) du \mid \; \, h \in \mathcal{S}(H)\right\}. \]

\end{definition}

\hop $H$ is said to be \emph{measurable} if its graph $\{(t,x) \mid \; \, x \in H(t)\}$ is measurable and \emph{integrally bounded} if there exists an integrable function $h : [0,T] \rightarrow \R_+$ such that

\[\sup_{x \in H(t)} \|x\| \leq h(t), \; \, \forall t \in [0,T].\]

\hop Let $h \in \mathcal{S}(H)$. We call $\Psi_h$ the map defined by $\tau \in [0,T]  \;  \mapsto \int_{0}^\tau h(u) du$.

\hop The following theorem is due to Aumann \cite{Aum65}

\begin{theoreme}\label{thmaumann} Let $H$ be a set-valued map on $[0,T]$ with nonempty images. Then

\begitem

\item[$*$] $\int_{[0,T]}  H(u) du$ is convex,

\item[$*$] If $H$ is measurable and integrally bounded then $\int_{[0,T]}  H(u) du$ is nonempty.

\item[$*$] If $H$ has closed images then $\int_{[0,T]} H(u) du$ is compact.

\item[$*$] If $(H_k)_k$ is a sequence of uniformly integrally bounded set-valued functions then

\[\limsup_k \int_{[0,T]} H_k(u) du \subset \int_{[0,T]} \limsup_k H_k (u) du,\]

\hop where $x \in \limsup_k A_k$ if and only if every neighborhood of $x$ intersects infinitely many $A_k$.

\end{itemize}

\end{theoreme}
\zun

\hop The next proposition is not a direct consequence of these results. However, the proofs of the third and fourth points can be adapted to derive it:

\begin{proposition} \label{convsimple} Let $(H_n)$ be a sequence of set-valued maps from $[0,T]$ to $\mathbb{R}^m$ uniformly integrally bounded  and $H$  be  a set-valued map with non empty images. We assume that, for any $\tau \in [0,T]$, $\limsup_n H_n(\tau) \subset H(\tau)$. For any $n \in \mathbb{N}$, let $h_n \in \mathcal{S}(H_n)$. Then,
\zun

\begitem

\item[$(i)$] If, for any $u$,  $H(u)$ is convex, there exists $h \in \mathcal{S}(H)$  such that $h_n$ converges weakly in $L^1([0,T])$ to $h$ (up to a subsequence). In particular, $\Psi_{h_n}$  converges simply to $\Psi_h$.
\item[$(ii)$] Without the convexity assumption, there exists a function $h$ on $[0,T]$ with the property that, for any $u \in [0,T]$, $h(u) \in co(H(u))$, the \emph{convex hull} of $H(u)$ (i.e. the smallest convex set containing $A$) and such that $h_n$ converges weakly in $L^1([0,T])$ to $h$. In particular, $\Psi_{h_n}$ converges simply to $\Psi_h$.


\end{itemize}

\end{proposition}
\zun

\hop {\bfseries Proof.} Since the sequence $(H_n)$ is uniformly integrally bounded, the $h_n$ are all bounded by an  integrable function $g : [0,T] \rightarrow \R_+$. Then there is a subsequence of $(h_n)$ with a weak limit $h \in L^1([0,T])$ (See Dunford and Schwartz \cite{DunSch58}, Theorem IV.8.9).We may assume without loss of generality that $(h_n)$ actually converges weakly to $h$. We now prove  that $h$ belongs to the set $\mathcal{S}(H)$.

\hop For $A\subset L^1([0,T])$, we call $\overline{co}(A)$ the smallest closed (for the $L^1$ norm) convex set containing $A$.

\hop Recall that, by Mazur's theorem, a convex subset of  $L^1([0,T])$ is closed if and only if it is weakly closed. Consequently, let $k \in \mathbb{N}$. the set $\overline{co}((h_n)_{n \geq k})$ is closed and convex and therefore weakly closed. Hence it contains $h$, which belongs to the weakly closed convex hull of $(h_n)_{n \geq k}$. hence,

\[h \in \overline{co}((h_n)_{n \geq k}) = \overline{co((h_n)_{n \geq k})},\]

\hop which means that there exists  $g_k \in co((h_n)_{n \geq k})$ such that $\|h - g_k\|_{L^1} \leq 1/k$. Finally, $(g_k)$ converges to $h$ in $L^1$ and we may assume without loss of generality that $(g_k)_k$ converges to $h$ almost everywhere on $[0,T]$. 

\hop From Caratheodory's theorem, the convex hull of a set $A$ is the set of all barycenters of families of $m+1$ elements of $A$.  Consequently, for any $k$ and $u \in [0,T]$, since $(h_n(u))_{n \geq k}$ is a set of points in $\mathbb{R}^m$, we have

\[g_k(u) = \sum_{j=0}^m \lambda_k^j(u) e_k^j(u),\]

\hop where $ \lambda_k^j(u) \geq 0$, $\sum_{j=0}^m \lambda_k^j(u) = 1$ and $e^j_k(u) \in \left\{h_n(u) \mid \; \, n \geq k \right\} \subset B(0,g(u))$. By compactness, we may assume that, for any $j$, $(e_k^j(u))_k$ converges to some $e^j(u)$ and $(\lambda_k^j(u))_k)$ converges to $\lambda^j(u)$, such that $\lambda_j(u) \geq 0$ and $\sum_{j=0}^m \lambda^j(u) = 1$. Finally,

\[h(u) = \lim_k g_k(u) = \sum_{j=0}^m \lambda^j(u) e^j(u).\]

\hop For any $j$, since $e^j(u)$ belongs to the limit set of the sequence $(h_n(u))_n$, it belongs to $H(u)$. Hence, $h(u)$ belongs to $co(H(u))$ and, if  $H(u)$ is convex, $h(u) \in H(u)$. The proof is complete. $\; \; \blacksquare$
\ztrois



\begin{lemme}\label{distancefixgen} Let $(K,d)$ be a compact metric space and $(\Lambda_n)_n$, $\Lambda : K \rightrightarrows K$ be set-valued maps such that $\Lambda$ is standard and $Fix(\Lambda) = \{x \in K \ | \ x \in \Lambda(x) \} \neq \emptyset$. Assume that, for any $x \in K$, 
\begitem
\item[$*$] $\Lambda_{n+1}(x) \subset \Lambda_n (x)$, 
\item[$*$] $\lim_n x_n = x \Rightarrow  \limsup_n \Lambda_n(x_n) \subset \Lambda(x).$
\end{itemize}

\hop Then, for all $\delta>0$, there exist $\epsilon>0$ and $n_0 \in \mathbb{N}$  such that for all $n> n_0$
$$d(x, \Lambda_{n}(x))\le \epsilon \ \Rightarrow  d(x, Fix (\Lambda)) \le \delta.$$

\end{lemme}
\zun

\hop {\bfseries Proof.}  First notice that, since $Fix(\Lambda)$ is non empty and $\Lambda_{n+1}(x) \subset \Lambda_n (x)$, there exists  $x \in K$ such that $x\in \Lambda_n(x)$ for all $n$.  Assume  that there are $\delta>0$, a strictly increasing sequence of integers $(n_k)_{k \geq 1}$ and a sequence $(x_k)_{k\ge 1}$ in $K$ such that
\[d(x_k, \Lambda_{n_k}(x_k)) \leq \frac{1}{k} \, \mbox{ and } \; \, d(x_k,Fix(\Lambda)) > \delta.\]

\hop Then, there exists a sequence $(y_k)_{k \geq 1}$ such that $y_k \in \Lambda_{n_k}(x_k)$ and
$$d(x_k, y_k) \le \frac{1}{k} \ \text{and }  d(x_k, Fix (\Lambda))> \delta.$$

\hop Without loss of generality, we may assume that $x_k \rightarrow x \in K$ and $y_k \rightarrow y$. Consequently, $d(x, Fix(\Lambda)) \geq  \delta> 0$ and $x=y$. On the other hand, 
\[y \in \limsup_k \Lambda_{n_k}(x_k) \subset \limsup_k \Lambda_k(x_k)  \subset \Lambda(x),\]

\hop which means that $x \in \Lambda(x)$, a contradiction. $\; \; \blacksquare$ 
\zdeux

\begin{Remarque} If the $\Lambda_n$ are closed, then it is sufficient to assume that $\cap_n \Lambda_n(x) = \Lambda(x), \; \forall x$. Indeed,  by monotonicity and the fact that $\Lambda_{n_k}$ is closed,  $y \in  \Lambda_{n_k}(x), \; \, \forall k$. 

\end{Remarque}

\hop From now, we consider a bounded standard set-valued map $F : \mathbb{R}^m \rightrightarrows \mathbb{R}^m$ and $T>0$.  Let $\delta$ be a positive real number. Then $F^{\delta}$ is the set-valued map defined by (\ref{fdelta}) and $\Lambda^{\delta}$ is the set-valued map defined by (\ref{lamdadelta}), extended to $\mathcal{C} \left([0,T], \mathbb{R}^m \right)$.
Note that, with our current notations,  the definition of the set-valued map $\Lambda^{\delta}$ can be written, for $\delta \geq 0$:
\[\Lambda^{\delta} : \mathcal{C} \left([0,T], \mathbb{R}^m \right) \rightrightarrows \mathcal{C} \left([0,T], \mathbb{R}^m \right), \; \;  \mathbf{z} \mapsto \left\{\mathbf{z}(0) + \Psi_h \mid h \in \mathcal{S}(F^{\delta}(\mathbf{z})) \right\}.\]

\begin{proposition} $\Lambda$ is a closed set-valued map with non empty images.
\end{proposition}
\zun

\hop {\bfseries Proof.} First, $\Lambda$ has non empty images since $\int_0^T F(\mathbf{z})$ is non empty, for any $\mathbf{z}$ (see Theorem \ref{thmaumann}). 

\hop Let $(\mathbf{z_n})_n$ be a sequence of $\mathcal{C}([0,T], \mathbb{R}^m)$, which converges to some $\mathbf{z}$ in $\mathcal{C}([0,T], \mathbb{R}^m)$ and let $(\mathbf{y_n})_n$ be a sequence converging to $\mathbf{y}$ such that, for all $n \in \mathbb{N}$, $\mathbf{y}_n \in \Lambda(\mathbf{z_n})$. This implies that, 
\[\forall n \in \mathbb{N}, \; \exists h_n \in \mathcal{S}(F(\mathbf{z_n})) \; \mbox{ such that } \; \, \mathbf{y_n}(\tau) = \mathbf{z_n}(0) + \Psi_{h_n}(\tau).\]

\hop We call $H_n := F(\mathbf{z_n})$ and $H := F (\mathbf{z})$. By assumptions we made on $F$, $H_n$ and H have compact, convex and nonempty values. For $\tau \in [0,T]$, $\mathbf{z_n}(\tau)$ converges to $\mathbf{z}(\tau)$ and, since the graph of $F$ is closed, 
\[\limsup_n H_n(\tau) = \limsup_n F (\mathbf{z_n}(\tau)) \subset F(\mathbf{z}(\tau)) = H(\tau).\]

\hop The assumptions of Proposition \ref{convsimple}  are satisfied. Hence, there exists $h \in \mathcal{S}(F(\mathbf{z}))$ such that $\mathbf{y} = \mathbf{z}(0) + \Psi_h$. $\; \; \blacksquare$
\zdeux

\begin{lemme} For any $\delta > 0$, $F^{\delta}$ is a closed set-valued map with non empty images.
\end{lemme}
\zun

\hop {\bfseries Proof.}  Let $x \in \mathbb{R}^m$. $F(x)$ is contained in $F^{\delta}$. Hence, it is not empty. Let $x_n \rightarrow x$ and  $(y_n)_n$ be a sequence of $F^{\delta}(x_n)$ converging to some $y$. Then there exists a sequence $(z_n)_n$ such that
\[d(z_n,x_n) < \delta \; \mbox{ and } \; \, d(y_n,F(z_n)) < \delta. \]

\hop hence there exists a sequence $(\alpha_n)_n$ such that $\alpha_n \in F (z_n)$ and $d(y_n,\alpha_n) < \delta$. Without loss of generality, we may assume that $z_n \rightarrow z$ and $\alpha_n \rightarrow \alpha$. By closeness of the graph of $F$, we obtain
\[d(z,x) < \delta, \; \, \alpha \in F(z) \, \mbox{ and } \; d(y,\alpha) < \delta,\]

\hop and $F^{\delta}$ is closed (and, in particular, has closed images). $\; \; \blacksquare$
\zun

\begin{Remarque} Note that the images are, a priori, not convex.

\end{Remarque}

\begin{lemme}\label{limsupfdelta} Let $(x_n)_n$ be a sequence of $\mathbb{R}^m$, converging to $x$ and $(\delta_n)_n$ be a positive, vanishing sequence. Then we have
\[\limsup_{n \rightarrow + \infty} F^{\delta_n}(x_n) \subset F(x).\]

\end{lemme}

\hop {\bfseries Proof.} Let $y \in \limsup_{n} F^{\delta_n}(x_n)$. By definition, there exists a sequence $(y_n)$ which converges to $y$ and such that $y_n \in F^{\delta_n}(x_n)$ (actually it is a subsequence but there is no loss of generality to keep the initial indexation). Hence there exists a sequence $(z_n)$ such that 
\[d(z_n,x_n) < \delta_n, \; \, d(y_n,F(z_n)) < \delta_n,\] 
which means that $d(x_n,\alpha_n) < \delta_n$ for some sequence $(\alpha_n)_n$ satisfying $\alpha_n \in F(z_n)$. Without loss of generality we may assume that $\alpha_n \rightarrow \alpha$ and $z_n \rightarrow z$. Hence we have  $y = \alpha \in F(z) = F(x)$.  $\; \; \blacksquare$
\zdeux

\begin{corollaire}\label{limsuplambdadelta} Let $(\mathbf{z_n})_n$ be a sequence converging to $\mathbf{z}$ in $\mathcal{C} \left([0,T], \mathbb{R}^m \right)$. Then we have, for any positive, vanishing sequence $(\delta_n)_n$,
\[\limsup_{n \rightarrow +\infty} \Lambda^{\delta_n}(\mathbf{z_n}) \subset \Lambda(\mathbf{z}).\]

\end{corollaire}

\hop {\bfseries Proof.} Let $\mathbf{y} \in \limsup_{n} \Lambda^{\delta_n}(\mathbf{z_n})$. This means that there exists a sequence $\mathbf{y_n} \in \Lambda^{\delta_n}(\mathbf{z_n})$ which converges to $\mathbf{y}$. Hence, for all $n \in \mathbb{N}$, there exists $h_n \in \mathcal{S}(F^{\delta_n}(\mathbf{z_n}))$ such that
\[\forall \tau \in [0,T], \; \mathbf{y_n}(\tau) = \mathbf{z_n}(0) + \int_0^{\tau} h_n(u) du\]

\hop By Corollary \ref{convsimple} and Lemma \ref{limsupfdelta}, there exists a function $h$ on $[0,T]$ such that 
\[\int_0^{\tau} h_n(u) du \rightarrow_n \int_0^{\tau} h(u) du, \; \, \forall \tau \in [0,T]\]

\hop and $h \in \mathcal{S}(F(\mathbf{z}))$, which completes the proof. $\; \; \blacksquare$
\zdeux

\hop 
\zun

\begin{corollaire} \label{distancefix} Let $\Lambda_n := \Lambda^{\delta_n}$. Suppose there exists a compact $K\subset \mathcal{C} \left([0,T], \mathbb{R}^m \right)$ such that $\Lambda_n : K \rightrightarrows K$ for all $n$. Then, for all $\delta>0$, there exists $\epsilon>0$ and $n_0 \in \mathbb{N}$ such that, for all $n>n_0$,
$$d(\mathbf{z}, \Lambda_{n}(\mathbf{z}))\le \epsilon \ \Rightarrow  d(\mathbf{z}, Fix (\Lambda)) \le \delta.$$

\end{corollaire}
\zun

\hop {\bfseries Proof.} This result follows from Lemma \ref{distancefixgen} and Corollary \ref{limsuplambdadelta}. $\; \; \blacksquare$
\zcinq

\hop {\bfseries Acknowledgements:} the authors would like to thank Michel Bena\"im for useful advices and discussions. They are also grateful to three anonymous referees as well as Eilon Solan for interesting remarks, in particular for suggesting the addition of Corollary \ref{potential}.

\bibliographystyle{plain}

\bibliography{biblio-jeuxdyn}

\end{document}